\documentstyle[12pt]{article}
\newcommand{\const}{\mathop{\rm const}\limits}

\newcommand{\supp}{\mathop{\rm supp}\limits}

\newcommand{\card}{\mathop{\rm card}\limits}

\newcommand{\Var}{\mathop{\rm Var}\limits}

\newcommand{\Law}{\mathop{\rm Law}\limits}

\newcommand{\grad}{\mathop{\rm grad}\limits}

\begin{document}

\begin{center}

{\bf  Non-asymptotical sharp exponential estimates for \\

\vspace{3mm}

maximum distribution of discontinuous random fields. } \\

\vspace{5mm}

 {\sc Ostrovsky E., Sirota L.}\\

\vspace{3mm}

 \ Department of Mathematics and Statistics, Bar-Ilan University, \\
59200, Ramat Gan, Israel.\\
e-mail: \ eugostrovsky@list.ru \\

\vspace{3mm}

 \ Department of Mathematics and Statistics, Bar-Ilan University,\\
59200, Ramat Gan, Israel.\\
e-mail: \ sirota3@bezeqint.net \\

\vspace{3mm}

 {\sc Abstract.}

\end{center}

 \  We offer in this paper the non-asymptotical bilateral sharp exponential estimates for tail of
maximum distribution of {\it discontinuous} random fields.\par

 \ Our consideration based on the theory of Prokhorov-Skorokhod spaces of random fields and  on the
 theory of multivariate Banach spaces of random variables with exponential decreasing tails of distributions. \par

 \vspace{3mm}

 {\it Key words and phrases:}  Random variable and random vector (r.v.), centered (mean zero) r.v. and  random processes
 r.p., moment generating function, rearrangement invariant Banach space on vector random variables, ordinary and exponential moments,
tail of distribution, Young-Orlicz function, norm, Chernov's estimate, upper and lower non-asymptotical exponential
estimates, Kramer's condition.

\vspace{4mm}

{\it Mathematics Subject Classification (2000):} primary 60G17; \ secondary
 60E07; 60G70.\\

\vspace{5mm}

\section{ Introduction. Previous results.} \par

\vspace{4mm}

 \hspace{3mm} We present here for beginning some known facts from the theory of one-dimensional random variables
with exponential decreasing tails of distributions, see  \cite{Ostrovsky4},  \cite{Kozachenko1}, \cite{Ostrovsky1}, chapters 1,2. \par

 \ Especially  we menton the authors preprint \cite{Ostrovsky5}; we offer in comparison with existing there results a more fine approach. \par

\vspace{3mm}
 \  Let $ (\Omega,F,{\bf P} ) $ be a probability space, $ \Omega = \{\omega\}. $ \par

\vspace{3mm}

 \ Let also $ \phi = \phi(\lambda), \lambda \in (-\lambda_0, \lambda_0), \ \lambda_0 =
\const \in (0, \infty] $ be certain even strong convex which takes positive values for positive arguments twice continuous
differentiable function, briefly: Young-Orlicz function,  such that
$$
 \phi(0) = \phi'(0) = 0, \ \phi^{''}(0) > 0, \ \lim_{\lambda \to \lambda_0} \phi(\lambda)/\lambda = \infty. \eqno(1.1)
$$

 \ For instance: $  \phi(\lambda) = 0.5 \lambda^2, \ \lambda_0 = \infty. $ \par

  \ We denote the set of all these  Young-Orlicz  function as $ \Phi; \ \Phi = \{ \phi(\cdot) \}. $ \par

 \ We say by definition that the {\it centered} random variable (r.v) $ \xi = \xi(\omega) $
belongs to the space $ B(\phi), $ if there exists some non-negative constant
$ \tau \ge 0 $ such that

$$
\forall \lambda \in (-\lambda_0, \lambda_0) \ \Rightarrow
\max_{\pm} {\bf E} \exp(\pm \lambda \xi) \le \exp[ \phi(\lambda \ \tau) ]. \eqno(1.2)
$$

 \ The minimal non-negative value $ \tau $ satisfying (1.2) for all the values $  \lambda \in (-\lambda_0, \lambda_0), $
is named a $ B(\phi) \ $ norm of the variable $ \xi, $ write

$$
||\xi||B(\phi)  \stackrel{def}{=}
$$

 $$
 \inf \{ \tau, \ \tau > 0: \ \forall \lambda:  \ |\lambda| < \lambda_0 \ \Rightarrow
  \max_{\pm}{\bf E}\exp( \pm \lambda \xi) \le \exp(\phi(\lambda \ \tau)) \}. \eqno(1.3)
 $$

 \ These spaces are very convenient for the investigation of the r.v. having a
exponential decreasing tail of distribution, for instance, for investigation of the limit theorem,
the exponential bounds of distribution for sums of random variables,
non-asymptotical properties, problem of continuous and weak compactness of random fields,
study of Central Limit Theorem in the Banach space etc.\par

  The space $ B(\phi) $ with respect to the norm $ || \cdot ||B(\phi) $ and
ordinary algebraic operations is a rearrangement invariant Banach space which is isomorphic to the subspace
consisting on all the centered variables of Orlicz's space $ (\Omega,F,{\bf P}), N(\cdot) $
with $ N \ - $ function

$$
N(u) = \exp(\phi^*(u)) - 1, \ \phi^*(u) = \sup_{\lambda} (\lambda u -
\phi(\lambda)).
$$
 \ The transform $ \phi \to \phi^* $ is called Young-Fenchel transform. The proof of considered
assertion used the properties of saddle-point method and theorem of Fenchel-Moraux:
$$
\phi^{**} = \phi.
$$

 \ The next facts about the $ B(\phi) $ spaces are proved in \cite{Kozachenko1}, \cite{Ostrovsky1}, p. 19-40:

$$
 \xi \in B(\phi) \Leftrightarrow {\bf E } \xi = 0, \ {\bf and} \ \exists C = \const > 0,
$$

$$
U(\xi,x) \le \exp(-\phi^*(Cx)), x \ge 0, \eqno(1.4)
$$
where $ U(\xi,x)$ denotes in this section the {\it one-dimensional tail} of
distribution of the r.v. $ \xi: $

$$
U(\xi,x) = \max \left( {\bf P}(\xi > x), \ {\bf P}(\xi < - x) \right), \ x \ge 0,
$$
and this estimation is in general case asymptotically as $ x \to \infty  $ exact. \par

 \ Here and further $ C, C_j, C(i) $ will denote the non-essentially positive
finite "constructive" constants.\par
 More exactly, if $ \lambda_0 = \infty, $ then the following implication holds:

$$
\lim_{\lambda \to \infty} \phi^{-1}(\log {\bf E} \exp(\lambda \xi))/\lambda =
K \in (0, \infty)
$$
if and only if

$$
\lim_{x \to \infty} (\phi^*)^{-1}( |\log U(\xi,x)| )/x = 1/K.
$$

 \ Hereafter $ f^{-1}(\cdot) $ denotes the inverse function to the
function $ f $ on the left-side half-line $ (C, \infty). $ \par

  \ Let $  F =  \{  \xi(t) \}, \ t \in T, \ T  $ is an arbitrary set, be the family of somehow
dependent mean zero random variables. The function $  \phi(\cdot) $ may be "constructive" introduced by the formula

$$
\phi(\lambda) = \phi_F(\lambda) \stackrel{def}{=} \max_{\pm} \log \sup_{t \in T}
 {\bf E} \exp(  \pm \lambda \xi(t)), \eqno(1.5)
$$
 if obviously the family $  F  $ of the centered r.v. $ \{ \xi(t), \ t \in T \} $ satisfies the  so - called
{\it uniform } Kramer's condition:
$$
\exists \mu \in (0, \infty), \ \sup_{t \in T} U(\xi(t), \ x) \le \exp(-\mu \ x),
\ x \ge 0.
$$
 In this case, i.e. in the case the choice the function $ \phi(\cdot) $ by the
formula (1.5), we will call the function $ \phi(\lambda) = \phi_0(\lambda) $
a {\it natural } function, and correspondingly the function

$$
\lambda \to {\bf E} e^{\lambda \xi}
$$
is named often as a moment generating function for the r.v.  $ \xi, $ if of course there exists  in some non-trivial
neighborhood of origin. \par

\vspace{3mm}

 \ The theory of multidimensional $  B(\phi) $ spaces is represented in the recent article
 \cite{Ostrovsky6}; it is quite analogous to the explained one. \par

 \ In detail,  denote by $  \epsilon  = \vec{\epsilon} = \{ \epsilon(1),   \epsilon(2), \ldots,   \epsilon(d) \} $ the non - random
$  d \ - $ dimensional  numerical vector, $  d = 2,3,\ldots, $ whose components take the values $  \pm 1 $ only.
 Set in particular $  \vec{1} = (1,1,\ldots,1) \in R^d_+. $ \par

 \ Denote also by $  \Theta = \Theta(d) = \{ \ \vec{\epsilon} \ \} $ {\it collection} of all such a vectors.  Note that
 $  \card \Theta = 2^d $  and $ \vec{1} \in \Theta.$  \par

 \ Another notations. For $ \vec{\epsilon} \in \Theta(d) $  and vector $ \vec{x} $ we introduce the coordinatewise product
as a $ d  \ - $ dimensional vector of the form

$$
\vec{\epsilon} \otimes \vec{x} \stackrel{def}{=} \{ \epsilon(1) \ x(1), \ \epsilon(2) \ x(2), \ \ldots,  \epsilon(d) \ x(d)  \}.
$$

\vspace{3mm}

 \ {\bf Definition 1.1.}\par

\vspace{3mm}

 \ Let $ \xi = \vec{\xi} = (\xi(1), \xi(2), \ldots, \xi(d) )  $ be a centered random vector such that each its component
$  \xi(j) $  satisfies the Kramer's condition. The {\it natural function} $ \phi_{\xi}= \phi_{\xi}(\lambda), \
 \lambda = \vec{\lambda} = (\lambda(1), \lambda(2), \ldots, \lambda(d))  \in R^d  $
for the random vector $  \xi  $ is defined as follows:

$$
\exp \{\phi_{\xi}(\lambda)\} \stackrel{def}{=} \max_{\vec{\epsilon}} \
{\bf E} \exp \left\{ \sum_{j=1}^d  \epsilon(j) \lambda(j) \xi(j) \right\}=
$$

$$
\max_{\vec{\epsilon} \in \Theta} \
{\bf E} \exp \{\epsilon(1) \lambda(1) \xi(1) + \epsilon(2) \lambda(2)\xi(2)+ \ldots +   \epsilon(d) \lambda(d)\xi(d) \}, \eqno(1.6)
$$
 where $  "\max" $ is calculated over all the combinations of signs $ \epsilon(j) =  \pm 1. $\par

\vspace{3mm}

 \ {\bf Definition 1.2.}\par

\vspace{3mm}

 \ The {\it tail function } for the random vector $ \vec{\xi} \ $
   $ U(\vec{\xi}, \vec{x}), \ \vec{x} = (x(1), x(2), \ldots, x(d)),  $ where all the coordinates $ x(j)  $
  of the deterministic vector $ \vec{x} $ are non-negative, is defined as follows.

$$
U(\vec{\xi}, \vec{x}) \stackrel{def}{=} \max_{\vec{\epsilon} \in \Theta}
{\bf P} \left( \cap_{j=1}^d  \{  \epsilon(j) \xi(j) > x(j) \}  \right) =
$$

$$
 \max_{ \vec{\epsilon} \in \Theta}
{\bf P}(\epsilon(1) \xi(1) > x(1), \ \epsilon(2) \xi(2) > x(2), \ \ldots, \ \epsilon(d) \xi(d) > x(d) ), \eqno(1.7)
$$
 where as before $  "\max" $ is calculated over all the combinations of signs $ \epsilon(j) =  \pm 1. $\par

 \ We illustrate this notion in the case $  d = 2.  $ Let $ \vec{\xi} = (\xi(1), \ \xi(2))  $ be a two-dimensional random
vector and let $ x, y  $  be non-negative numbers. Then

$$
U( (\xi(1), \xi(2)), \ (x,y) ) =
$$

$$
\max [ {\bf P} (\xi(1) > x, \ \xi(2) > y), \ {\bf P} (\xi(1) > x, \ \xi(2) < - y),
$$

$$
{\bf P} (\xi(1) < - x, \ \xi(2) > y), \ {\bf P} (\xi(1) < - x, \ \xi(2) < - y) ].
$$

\vspace{3mm}

 \ {\bf Definition 1.3.}\par

\vspace{3mm}

 \ Let $  h = h(x), \ x \in R^d  $ be some non - negative real valued function, which is finite on some non-empty
neighborhood of origin. We denote as ordinary

$$
\supp h = \{x, \ h(x) < \infty   \}.
$$

 \ The Young-Fenchel, or Legendre transform $  h^*(y), \ y \in R^d  $ is defined likewise the one-dimensional case

$$
h^*(y) \stackrel{def}{=} \sup_{x \in \supp h} ( (x,y) - h(x)). \eqno(1.8)
$$

 \ Herewith $ (x,y)  $ denotes the scalar product of the vectors $ x,y: \ (x,y) = \sum_j x(j)y(j);
  \ |x| = \sqrt{(x,x)}. $\par

 \ Obviously, if the set $ \supp h $ is central symmetric, then the function $ h^*(y) $ is even.\par

\vspace{3mm}

 \ {\bf Definition 1.4.}\par

\vspace{3mm}

 \ Recall, see \cite{Krasnosel'skii1}, \cite{Rao1}, \cite{Rao2} that the function $  x \to g(x), \ x \in R^d, \ g(x) \in R^1_+  $ is named
 multivariate Young, or Young-Orlicz function, if it is even,  $ d \ - $ times  continuous differentiable, convex, non-negative,
 finite on the whole space $  R^d,$ and such that

$$
g(x) = 0 \ \Leftrightarrow x = 0; \hspace{4mm} \frac{\partial g}{\partial x}/(\vec{x} = 0) = 0,
$$

$$
\det \frac{\partial^2 g}{\partial x^2} /(\vec{x} = 0) > 0. \eqno(1.8)
$$

 We explain in detail:

$$
 \frac{\partial g}{\partial x} = \left\{  \frac{ \partial g}{ \partial x_j} \right\} = \grad g,  \hspace{5mm}
\frac{\partial^2 g}{\partial x^2} = \left \{ \frac{\partial^2 g}{\partial x_k \partial x_l}  \right\}.
\eqno(1.9)
$$

 \ We assume  finally

$$
\lim_{|x| \to \infty} \frac{\partial^d g }{\prod_{k=1}^d  \partial x_k} = \infty. \eqno(1.10)
$$

 \ We will denote the set of all such a functions by $  Y = Y(R^d) $ and denote also by $ D  $ introduced before matrix

$$
D = D_g := \frac{1}{2} \left\{ \frac{\partial^2 g(0)}{\partial x_k \partial x_l} \right\}.
$$

 \ Evidently, the matrix $ D = D_g $ is non-negative definite,  write $  D = D_g \ge \ge 0.  $ \par

\vspace{3mm}

 \ {\bf Definition 1.5.}\par

\vspace{3mm}

  \ Let $  V, \ V \subset R^d $ be open convex  central symmetric:  $ \forall x  \in V \ \Rightarrow -x \in V  $
 subset of whole space $  R^d $ containing some non-empty neighborhood of origin. We will denote the  collection
 of all such a sets  by $  S = S(R^d).  $\par

   \ We will denote by $  Y(V), \ V \in S(R^d) $ the set of all such a functions from the definition  2.4 which are defined only
on the set $  V. $ For instance, they are even, twice  continuous differentiable, convex, non-negative, is equal to zero
only at the origin  and $  \det D_g > 0,  $

$$
\lim_{x \to \partial V -0} \frac{\partial^d g }{\prod_{k=1}^d  \partial x_k} = \infty. \eqno(1.11)
$$

 \ Notation: $  V = \supp g.  $ \par

\vspace{3mm}

 \ {\bf  Definition 1.6. } \par

\vspace{3mm}

 \ Let the set $ V  $ be from the set
  $ S(R^d): V \in S(R^d)  $ and let the Young function  $  \phi(\cdot)  $ be from the set $  Y(V): \supp \phi = V. $\par

 \ We will say by definition likewise the one-dimensional case  that the {\it centered} random vector (r.v)
 $ \xi = \xi(\omega) = \vec{\xi} = (\xi(1), \xi(2), \ldots, \xi(d)) $ with values in the space $  R^d $ belongs to the
 space $ B_V(\phi), $  write $ \vec{\xi}\in  B_V(\phi),  $  if there exists certain non-negative constant $ \tau \ge 0 $ such that

$$
\forall \lambda \in V \ \Rightarrow
\max_{\vec{\epsilon}} {\bf E} \exp \left( \sum_{j=1}^d \epsilon(j) \lambda(j) \xi(j) \right) \le
\exp[ \phi(\lambda \cdot \tau) ]. \eqno(1.12)
$$

 \ The minimal value $ \tau $ satisfying (3.1) for all the values $  \lambda \in V, $
is named by definition as a $ B_V(\phi) \ $ norm of the vector $ \xi, $ write

$$
||\xi||B_V(\phi)  \stackrel{def}{=}
$$

 $$
 \inf \left\{ \tau, \ \tau > 0: \ \forall \lambda:  \ \lambda \in V \ \Rightarrow
 \max_{\vec{\epsilon}}{\bf E}\exp \left( \sum_{j=1}^d \epsilon(j) \lambda(j) \xi(j) \right) \le
 \exp(\phi(\lambda \cdot \tau)) \right\}. \eqno(1.13)
 $$

\vspace{3mm}

 \ The space $ \ B_V(\phi) \ $  relative introduced here norm $ ||\xi||B_V(\phi)  $ and ordinary algebraic operation
 is also rearrangement invariant Banach space. \par

\vspace{3mm}

 \ For example, the {\it generating function } $  \phi_{\xi}(\lambda) $ for these spaces may be picked by the
 following natural way:

$$
\exp[ \phi_{\xi}(\lambda ) ] \stackrel{def}{=}
\max_{\vec{\epsilon} \in \Theta} {\bf E} \exp \left( \sum_{j=1}^d \epsilon(j) \lambda(j) \xi(j) \right),
 \eqno(1.13a)
$$
if of course the random vector  $  \xi $ is centered and has an exponential tail of distribution. This imply  that
the natural function  $ \phi_{\xi}(\lambda ) $ is finite on some non-trivial central symmetrical neighborhood of origin,
or equivalently  the mean zero random vector  $  \xi $ satisfies the multivariate Kramer's condition. \par

 \ Obviously, for the natural function $ \phi_{\xi}(\lambda )  $

$$
||\xi||B(\phi_{\xi}) = 1.
$$

 \ It is easily to see that this  choice of the generating function $ \phi_{\xi} $ is optimal, but in the practical using
often this function can not  be calculated in explicit view,  but there is a possibility to estimate its. \par

 \vspace{3mm}

 \ We agree to  take in the case when   $  V = R^d \hspace{4mm}  \ B_{R^d}(\phi):= B(\phi).  $ \par

 \ Note that the expression for the norm $ ||\xi||B_V(\phi)  $ dependent aside from the function $  \phi  $ and the set $ V, $
only on the distribution $  \Law(\xi). $ Thus, this norm and correspondent space  $  B(\phi)  $ are rearrangement invariant (symmetrical)
in the terminology of the classical book  \cite{Bennet1}, see chapters 1,2. \par

\vspace{4mm}

 \ The following important facts about upper tail estimate is proved in the preprint \cite{Ostrovsky6}. \par

\vspace{4mm}

 \ Let $  \phi = \phi(\lambda), \ \lambda \in V \subset R^d $ be arbitrary non - negative real valued function,
which is finite on some  non - empty symmetrical neighborhood of origin.  Suppose for given centered $  d  \ - $ dimensional random vector $ \xi = \vec{\xi}  $

$$
{\bf E} e^{(\lambda,\xi)} \le e^{ \phi( \lambda) }, \ \lambda \in V.  \eqno(1.14)
$$

 \ On the other words, $ ||\xi||B(\phi) \le 1. $ \ Then for all non-negative vector $ x = \vec{x}  $ there holds

 $$
 U(\vec{\xi}, \vec{x}) \le \exp \left( - \phi^*(\vec{x}) \right) \ - \eqno(1.15)
 $$
the multidimensional generalization of Chernov's  inequality. \par
 \ Moreover, the last estimate is essentially non-improvable; there are some lower estimates for considered here
 multivariate  tail function in \cite{Ostrovsky6} . \par

\vspace{4mm}

 \ {\bf  Corollary 1.1. }   Let as above the function $  \phi(\cdot)  $ be from the Young-Orlicz set and
such that $   \phi(0) = 0. $ The centered non-zero random vector $  \xi  $ belongs to the space  $   B(\phi): $

$$
\exists C_1 \in (0,\infty), \ \forall \lambda \in R^d
 \Rightarrow  {\bf E} e^{(\lambda,\xi)} \le e^{ \phi( C_1 \cdot \lambda) }, \ \lambda \in R^d \eqno(1.16)
$$
if and only if

$$
\exists C_2 \in (0,\infty), \ \forall \ x \in R^d_+  \Rightarrow \  U(\vec{\xi}, \vec{x}) \le \exp \left( - \phi^*(\vec{x}/C_2) \right).
\eqno(1.17)
$$

 \ More precisely, the following implication holds: there is finite positive constant $  C_3 =  C_3(\phi)  $ such that
  for arbitrary non - zero centered r.v. $  \xi: \ ||\xi|| = ||\xi||B(\phi)  < \infty \ \Leftrightarrow $

$$
 \forall \lambda \in R^d
 \Rightarrow  {\bf E} e^{(\lambda,\xi)} \le e^{ \phi( ||\xi|| \cdot \lambda) }
$$
iff

$$
\exists C_3(\phi) \in (0,\infty) \ \forall \ x \in R^d_+  \Rightarrow \  U(\vec{\xi}, \vec{x}) \le
\exp \left( - \phi^*(\vec{x}/(C_3 / ||\xi||) \right). \eqno(1.18)
$$

\vspace{3mm}

\ {\bf  Corollary 1.2. } Assume the non-zero centered random vector $  \xi = (\xi(1), \xi(2), \ldots, \xi(d)) $ belongs
to the space $  B(\phi):  $

$$
{\bf E} e^{(\lambda,\xi)} \le e^{ \phi(||\xi|| \cdot \lambda) }, \ \phi \in Y(R^d),
\eqno(1.19)
$$
and let $  y  $ be arbitrary positive non-random number. Then $  \forall y  > 0 \ \Rightarrow  $

$$
{\bf P} \left( \min_{j = 1,2,\ldots,n} |\xi(j)| > y  \right) \le 2^d \cdot
\exp \left(- \phi^*(y/||\xi||,y/||\xi||, \ldots, y/||\xi||) \right).  \eqno(1.20)
$$

\vspace{4mm}

{\bf Example 1.1.} Let as before $  V = R^d  $ and $  \phi(\lambda) = \phi^{(B)}(\lambda)  = 0.5(B\lambda,\lambda),  $ where
$  B  $  is non-degenerate positive definite symmetrical matrix, in particular $  \det B > 0. $ It follows from theorem 6.1
that the (centered) random vector $ \xi $ is subgaussian relative the matrix $ B: $

$$
\forall \lambda \in R^d \ \Rightarrow {\bf E} e^{ (\lambda, \xi)} \le e^{0.5 (B \lambda, \lambda) ||\xi||^2 }.
$$
iff for some finite positive constant $  K = K(B,d) $ and for any  non - random positive vector $  x = \vec{x} $

$$
U(\xi,x) \le e^{- 0.5 \ \left( (B^{-1}x,x)/(K||\xi||^2) \right) }. \eqno(1.21)
$$

\vspace{3mm}

\section{Positive part of our exponential spaces.}

\vspace{3mm}

 \hspace{6mm} Let $ (\Omega,F,{\bf P} ) $ be again a probability space, $ \Omega = \{\omega\}. $ \par

 \vspace{3mm}

\ Let also as before $ \phi = \phi(\lambda), \lambda \in (-\lambda_0, \lambda_0), \ \lambda_0 =
\const \in (0, \infty] $ be certain even strong convex which takes positive values for positive arguments twice continuous
differentiable function, such that
$$
 \phi(0) = \phi'(0) = 0, \ \phi^{''}(0) > 0, \ \lim_{\lambda \to \lambda_0} \phi(\lambda)/\lambda = \infty.
$$
  \ We denote the set of all these  (Young-Orlicz)  functions as $ \Phi; \ \Phi = \{ \phi(\cdot) \}. $ \par

\vspace{3mm}

 \ {\bf Definition 2.1.} \par

 \ {\it We say by definition that the  random variable (r.v) $ \xi = \xi(\omega) $
belongs to the space $ B^+(\phi):  \ \xi \in B^+(\phi), $ if there exists some non-negative constant
$ \tau \ge 0 $ such that}

\vspace{3mm}

$$
\forall \lambda \in (0, \lambda_0) \ \Rightarrow
 {\bf E} \exp( \lambda \xi) \le \exp[ \phi(\lambda \ \tau) ]. \eqno(2.1)
$$

 \ The minimal non - negative value $ \tau $ satisfying this restriction for all the values $  \lambda \in (0, \lambda_0), $
is named a $ B^+(\phi) \ $ norm, more precisely, quasy-norm, of the variable $ \xi, $ write

$$
||\xi||B^+(\phi) = ||\xi||^+   \stackrel{def}{=}
$$

 $$
 \inf \{ \tau, \ \tau > 0: \ \forall \lambda: 0 \le  \ \lambda < \lambda_0 \ \Rightarrow
  {\bf E}\exp( \lambda \xi) \le \exp(\phi(\lambda \ \tau)) \}. \eqno(2.2)
 $$

 \ The theory of these semi-Banach space  is represented in the monograph \cite{Ostrovsky1},
chapter 1, sections 1.1; 1.2; chapter 3, section 3.16.\par
  \ Note tat in \cite{Ostrovsky1} it was considered only the case when $  \phi(\lambda) = 0.5 \cdot \lambda^2, \ \lambda \in R. $ \par

\vspace{3mm}

 \ Let us itemize some  main properties of these quasy-Banach spaces.

$$
{\bf 1.} \ \xi \in B^+(\phi) \ \Rightarrow {\bf E} \xi \le 0.
$$

$$
{\bf 2.} \ \xi \in B^+(\phi), \ \lambda \ge 0 \ \Rightarrow \lambda \xi \in B^+(\phi)
$$
and herewith

$$
||\lambda \ \xi||B^+(\phi) = \lambda || \ \xi||B^+(\phi).
$$

$$
{\bf 3.} \ \xi, \ \eta \in B^+(\phi) \Rightarrow  \xi + \eta \in B^+(\phi).
$$

$$
{\bf 4.} \ \xi, \ \eta \in B^+(\phi), \ \phi \in \Phi  \ \Rightarrow  \ ||\xi + \eta||^+  \le
||\xi||^+  + ||\eta||^+.
$$

  \ If in addition the function $ \lambda \to \phi( \sqrt{\lambda}), 0 \le \lambda < \lambda_0 $ is convex
 and both the r.v. being from the space $  B^+\phi, $ say $ \xi, \ \eta $ are independent, then

$$
{\bf 4a.} \ \xi, \ \eta \in B^+(\phi), \ \phi \in \Phi  \ \Rightarrow  \ ||\xi + \eta||^+  \le
\sqrt{(||\xi||^+)^2 +  (||\eta||^+)^2}.
$$

$$
{\bf 5.} \
\xi \in B^+(\phi) \ \Longleftrightarrow \exists C \in (0, \infty], \ V(\xi, x) \stackrel{def}{=}
{\bf P} (\xi > x) \le  \exp \left(-\phi^*(x/C) \right).
$$

 \ Hereafter $  \phi^*(\cdot) $ denotes the classical Legendre, or Young-Fenchel transform of the function $ \phi(\cdot): $

$$
\phi^*(z) \stackrel{def}{=} \sup_x ( (x,z) - \phi(x)).
$$

\vspace{4mm}

 \ {\bf An important example.}

\vspace{3mm}

 \ Let $  \phi = \phi(\lambda), \ \lambda = \vec{\lambda} \in V \subset R^d $ be certain non-negative real valued twice
continuous  differentiable Young-Orlicz function, which is finite on some-empty symmetrical relative the origin
neighborhood.  Suppose for given centered $  d  \ - $ dimensional random vector $ \xi = \vec{\xi} = \{\xi(1), \xi(2), \ldots, \xi(d) \} $ there
exists a non-negative constant $ \tau \ge 0 $ such that

$$
\forall \lambda \in V \ \Rightarrow
{\bf E} e^{(\lambda,\xi)} \le e^{ \phi( \tau \cdot \lambda) }.
$$
 The notation $   (\lambda,\xi)  $ denotes as usually  the scalar (inner) product of both the vectors $ \lambda $  and $ \xi.$\par

 \ The minimal non-negative value $ \tau $  satisfying the last inequality for all the admissible values $  \vec{\lambda} \in V,  $
is names as before $  B(\phi) \ - $ norm of the random vector $  \vec{\xi}; $  write: $  ||\xi||B(\phi) = \inf \tau. $ \par
 \ The r.v. belonging to this Banach rearrangement invariant space  $  B(\phi) $  having exponential decreasing tail of distribution. \par
 \ The detail investigation of these spaces may be found in the recent article \cite{Ostrovsky6}. \par

\vspace{3mm}

 \ Assume now the non-zero centered random vector $  \xi = (\xi(1), \xi(2), \ldots, \xi(d)) $ belongs
to the space $  B(\phi):  $

$$
{\bf E} e^{(\lambda,\xi)} \le e^{ \phi(||\xi|| \cdot \lambda) },
$$
and let $  y  $ be arbitrary positive non-random number. Then $  \forall y  > 0 \ \Rightarrow  $

$$
{\bf P} \left( \min_{j = 1,2,\ldots,d} \ |\xi(j)| > y  \right) \le 2^d \cdot
\exp \left(- \phi^*(y/||\xi||,y/||\xi||, \ldots, y/||\xi||) \right). \eqno(2.3)
$$

\vspace{4mm}

\section{ Exact exponential estimates for maximum distribution of discontinuous random fields. }

\vspace{3mm}

 \ {\bf A. Upper bounds.}\par

\vspace{3mm}

 \ We are passing to the main goal of this article, applying the results of the last section, especially the estimation (2.3).

 \vspace{3mm}

 \ Let $   T = [0,1]^m  $ be an ordinary $ m \ -  $ dimensional unit cube and let

$$
\xi = \xi(\vec{t}) = \xi( t(1), t(2), \ldots, t(m) ), \hspace{4mm} t(i) \in [0,1], \ t = \vec{t} \in T \eqno(3.0)
$$
be a separable numerical valued random process (field).  Put

$$
P(T,u) = {\bf P} ( \sup_{t \in T} \xi(t) > u). \eqno(3.1)
$$

 \ {\it Our goal is obtaining an exact  non-asymptotical  exponential  decreasing to zero as $  u \to \infty  $
estimates for } $  P(T,u). $  \par

\vspace{3mm}

 \ This problem can be considered as a classical; see for instance \cite{Buldygin3}, \cite{Fernique1}, \cite{Fernique2},
\cite{Fernique3}, \cite{Kahane1}, \cite{Ostrovsky1}, \cite{Pizier1}, \cite{Talagrand1} etc.\par
 \ Note that in the mentioned articles and books the considered random processes (and fields)  are continuous almost
everywhere.  We presume to investigate in this article the discontinuous r.f.  More precisely, it follows under imposed us
conditions that the r.f $  \xi(\vec{t}) $ belongs only the Prokhorov-Skorokhod  space $  D[0,1]^m,  $
\cite{Prokhorov2}, \cite{Skorokhod1}, \cite{Kolmogorov1}. See also \cite{Bickel1}, \cite{Billingsley1}, \cite{Billingsley2}, \cite{Gikhman1}. \par

\ The first well known application of considered estimates are in the non-parametric statistics, in the theory of empirical distribution
function, see \cite{Prokhorov2}, \cite{Skorokhod1}, \cite{Bickel1}, \cite{Neuhaus1}. The second one may be found  in the method Monte-Carlo
by computation of multiparameter integrals, \cite{Frolov1}, especially by  computation of  multiparameter integrals with discontinuous
integrand function  \cite{Grigorjeva1}. \par

\vspace{4mm}

 \ We need to introduce some another notations.
 Denote by $  \epsilon  = \vec{\epsilon} = \{ \epsilon(1),   \epsilon(2), \ldots,   \epsilon(d) \} $ the non-random
$  d \ - $ dimensional  numerical vector, $  d = 2,3,\ldots, $ whose components take the values $  \pm 1 $ only.  Set
$  \vec{1} = (1,1,\ldots,1) \in R^d_+. $ \par

 \ Denote by $  \Theta = \Theta(d) = \{ \ \vec{\epsilon} \ \} $ {\it collection} of all such a vectors.  Note that
 $  \card \Theta = 2^d $  and $ \vec{1} \in \Theta.$  \par

 \ For the value $ \vec{\epsilon} \in \Theta(d) $  and vector $ \vec{x} $ we introduce the coordinatewise product
as a $ d  \ - $ dimensional vector of the form

$$
\vec{\epsilon} \otimes \vec{x} \stackrel{def}{=} \{ \epsilon(1) \ x(1), \ \epsilon(2) \ x(2), \ \ldots,  \epsilon(d) \ x(d)  \}.
$$

\vspace{4mm}

 \ Let for beginning $  T = [0,1].  $  An arbitrary finite  subset $  T_N = \{  s(1), s(2), \ldots, s(N) \} $
  of the set $  T: \ T_N \subset T  $ is said to be  {\it net}  of $  T.  $ We require without loss of generality
$  0 \le s(1) < s(2) < s(3) < \ldots s(N) \le 1. $  Note that $  N = \card(T_N). $ \par

 \ Let $  s $ be arbitrary point from the set $ T = [0,1] $ and let $   T_N $  be any net of $  T. $ Define the following
 functions ("projections") of the arbitrary point $  s  $ from the set $  T:  $
 $  \theta^{+1}(s), \ \theta^{-1}(s)  = \theta^{+1}(T_N,s), \ \theta^{-1}(T_N,s):  $

$$
\theta^{+1}(s) = s(k) \in T_N, \ s(k) = \min_j \{s(j): \ s(j) \ge s  \};
$$

$$
\theta^{-1}(s) = s(l) \in T_N, \ s(l) = \max_i \{s(i): \ s(j) < s  \}.
$$

 \ One can $ \theta^{-1}(s) = 0  $ or $ \theta^{+1}(s) = 1. $ \par

 \ So, we define the function $  \theta^{\epsilon}(T_N,s)  $ for arbitrary net $  T_N, $  for the values $ s \in T $
and for $  \epsilon = \pm 1. $ \par
 \ Let us passing to the multidimensional case. Let $  \vec{t} = \{ t(1), t(2), \ldots, t(m) \} $ be an arbitrary point from the unit
cube $  T = [0,1]^m  $  and $  \vec{\epsilon} = \epsilon \in \Theta_m. $  Let also
$ \vec{T}_{\vec{N}} = \{ T_{N(j)} \} $ be a direct (Descartes) product of the $ \ m - $ tuples of the "individual"  $ T_{N(j)} $
  nets in the interval $ [0,1]: $

$$
\vec{T}_{\vec{N}} := \otimes_{j=1}^m T_{N(j)}.
$$

 \ We agree to take the following multi-value "projections"
 of arbitrary  point $ t = \vec{t} \in \vec{T}  $
 on the net $ \vec{T}_{\vec{N}} $

$$
\theta^{\vec{\epsilon} }( \vec{t} ) = \theta^{\vec{\epsilon} }(\vec{T}_{\vec{N}}, \vec{t} )
 \stackrel{def}{=} \{\theta^{\epsilon(1)} t(1), \  \theta^{\epsilon(2)} t(2), \ldots,
    \theta^{\epsilon(m)} t(m) \}, \eqno(3.2)
$$
and define the optimal "projection"  $ \theta (\vec{t}) = \theta(\vec{T}_{\vec{N}} , \vec{t})   $
of arbitrary  point $ t = \vec{t} \in \vec{T}  $ on the net $ \vec{T}_{\vec{N}} $
as follows:

$$
|| \left[ \xi(\vec{t}) - \xi \left(\theta(\vec{t}) \right) \ \right] \  ||^+B(\vec{\phi})=
 || \min_{\epsilon \in \Theta(m)} \left[ \xi(\vec{t}) - \xi \left(\theta^{\vec{\epsilon}}(\vec{t}) \right) \ \right] \  ||^+B(\vec{\phi}).
 \eqno(3.3)
$$

\ The possible non-uniqueness of the value $ \theta (\vec{t}) $ does not play a role; it is important only to choose
this variable deterministic. \par

\ The exponential estimations for the distribution of minimum for finite set of random variables is described in the second section. \par

 \ Note that the net $ \vec{T}_{\vec{N}}  $ on the whole set $  T = \vec{T} = [0,1]^m  $ contains

 $$
 K(\vec{T}_{\vec{N}}):= \card( \vec{T}_{\vec{N}}) = \prod_{j=1}^m N(j)
 $$
elements. \par

 \ Let us introduce the following important characteristic for the defined net $ \vec{T}_{\vec{N}} $

$$
\Delta( \vec{T}_{\vec{N}}) \stackrel{def}{=} \sup_{\vec{t} \in [0,1]^m} \
|| \min_{\epsilon \in \Theta(m)} \left[ \xi(\vec{t}) - \xi \left(\theta^{\vec{\epsilon}}(\vec{t}) \right) \ \right] \  ||^+B(\vec{\phi})=
$$

$$
 \sup_{\vec{t} \in [0,1]^m} \
|| \left[ \xi(\vec{t}) - \xi \left(\theta(\vec{t}) \right) \ \right] \  ||^+B(\vec{\phi}),  \eqno(3.4)
$$
and the "inverse" function

$$
M(\delta) = M(\phi, \delta) \stackrel{def}{=}
\inf \left\{ \card \left( \vec{T}_{\vec{N}} \right): \Delta \left( \vec{T}_{\vec{N}} \right) < \delta \right\}, \
\delta > 0.  \eqno(3.4a)
$$

 \ Denote also $ m(\delta) = m(\phi, \delta) = \ln  M(\delta) = \ln M(\phi, \delta), $

$$
g(p) = g(p,\phi):= (1 - p) \sum_{n=1}^{\infty}  p^{n-1} m(p^n) = (1 - p) \sum_{n=1}^{\infty}  p^{n-1} m(\phi, p^n).   \eqno(3.5)
$$

 \  This values are the direct extension of the classical definitions of a metric entropy and correspondingly entropy series,
which was used by an investigation of continuous random fields, see \cite{Fernique1}, \cite{Kozachenko1}, \cite{Ostrovsky1},
\cite{Pizier1}. \par

\vspace{4mm}

{\bf Theorem A.} Let $  \phi = \phi(\lambda) $ be Young-Orlicz  function, $  \phi(0) = \phi(0+) = 0,  $ and such that
for the  centered r.f. $ \xi(t), \ t \in [0,1]^m  $

$$
 \sup_{t \in T} {\bf E} e^{\lambda \xi(t)} \le e^{\phi(\lambda) }, \ \lambda \ge 0. \eqno(3.6)
$$

 \ It will be presumed the  finiteness of these function on some right-hand side non-empty neighborhood of zero. \par

 \ For instance, the function $ \phi = \phi_0(\lambda) $ may be picked by the natural way:

$$
e^{\phi_0(\lambda) } :=  \sup_{t \in T} {\bf E} e^{\lambda \xi(t)}, \ \lambda \ge 0. \eqno(3.6a)
$$

\ Let also

$$
 \exists p \in (0,1)  \ \Rightarrow g(p) < \infty. \eqno(3.7)
$$

\vspace{4mm}

\ We propose

$$
P(T,u) \le \inf_{p \in (0,1) } \exp \left[ g(p) - \phi^*(u(1-p)) \right], \ u > 0. \eqno(3.8)
$$

\vspace{3mm}

{\bf Proof.} \par
 \ Denote for brevity $  \zeta = \sup_t \xi(t) $ and
 for any random variable $ \ \eta \hspace{4mm}  ||\eta ||  = ||\eta ||B^+(\phi). $ \par

 \ One can suppose without loss of generality the existence of the point $  t_0 \in T $ such that $ ||\xi(t_0) || = 1.  $
 \ Let  $  p = \const \in (0,1); $ the exact value of $ p   $ will be clarified below. \par

 \ Let us introduce also the following non - decreasing (imbedded) sequence $ S(n), \ n = 0,1,2,\ldots  $  of the nets inside the set
$  T: \ S(0) := \{  t_0 \},  \ $ and the net $ S(n)  $ is such that $  \card(S(n)) =  M(p^n).  $ The described before {\it optimal }
projection of the point $  t \in T $ onto the net $  S(n) $ will be denoted  again by $ \theta_n t; \ \theta_n t \in S(n). $\par
 \ By virtue of separability of the r.f. $  \xi(\cdot) $

$$
\zeta = \sup_{t \in T} \xi(t) = \lim_{n \to \infty} \max_{t \in S(n)} \xi(t),
$$

$$
 \max_{t \in S(n)} \xi(t) \le  \max_{t \in S(n)} \left[ \xi(t) - \xi(\theta_{n-1}t) \right] + \max_{t \in S(n-1)} \xi(t), \ n \ge 1.
$$

 \ We conclude therefore denoting $  \eta_n = \max_{t \in S(n)} \left[ \xi(t) - \xi(\theta_{n-1}t) \right] $

$$
\sup_{t \in T} \xi(t) \le \sum_{n=1}^{\infty} \eta_n.  \eqno(3.9)
$$

 \ We apply the famous H\"older's inequality

$$
{\bf E }e^{ \lambda \ \sup_{t \in T} \xi(t) } \le \prod_{n=1}^{\infty} \left[ {\bf E} e^{\lambda \ r_n \ \eta_n} \right]^{1/r_n},
$$
 where

$$
\lambda > 0, \ r_n > 1, \ \sum_{n=1}^{\infty}  1/r_n = 1.
$$

 \ Let us estimate each factor.

$$
{\bf E}e^{\lambda \ \eta_n} = {\bf E} \max_{t \in S(n)} e^{\lambda (\xi(t) - \xi(\theta_{n-1} t ) ) } \le
\sum_{t \in S(n)} {\bf E} e^{\lambda (\xi(t) - \xi(\theta_{n-1} t) )} \le
$$

$$
\sum_{t \in S(n)} \ e^{ \phi \left(\lambda \ p^{n-1} \right) } = M(p^n) \  e^{ \phi \left(\lambda \ p^{n-1} \right) },
$$
following

$$
\ln {\bf E} e^{\lambda \ \zeta} \le \sum_{n=1}^{\infty}  \frac{m(p^n)}{r_n} + \sum_{n=1}^{\infty}  \frac{\phi(\lambda \ p^{n-1} r_n)}{r_n}.
$$

 \ Let us choose $  r_n = \frac{1}{(1-p) \ p^{n-1}}, $ then we obtain after some simple calculations

 $$
 \ln {\bf E} e^{\lambda \ \zeta} \le (1-p) \sum_{n=1}^{\infty} p^{n-1} \phi\left(\frac{\lambda}{1-p} \right) +
(1-p) \sum_{n=1}^{\infty}  p^{n-1} m(p^n) =
 $$

$$
\phi\left(\frac{\lambda}{1-p} \right) + g(p).  \eqno(3.10)
$$

 \ We conclude finally using Chernov's inequality

$$
P(T,u) \le \inf_{p \in (0,1) } \exp \left[ g(p) - \phi^*(u(1-p)) \right], \ u > 0.  \eqno(3.11)
$$

\vspace{4mm}

 \ Let us consider some examples. \par

\vspace{3mm}

 \ {\bf Example 1.} Define

$$
\pi(u) = \frac{1}{u \cdot \phi^{*'}(u)}
$$
for the sufficiently greatest values $ u, \ u \ge u_0, $ where $  u_0 = \const $ for which $ \pi(u_0) = p_0,  $
say for definiteness $  p_0 = 1/2.$  We obtain for these values $ u $

$$
P(T,u) \le \exp \left\{- \phi^*(u) + g(\pi(u)) + 1/2  \right\}.
$$
 Of course, the last estimate is meaningful only in the case when

$$
\overline{\lim}_{u \to \infty} \frac{g(\pi(u))}{\phi^*(u)} <1.
$$

\vspace{3mm}

\ {\bf Example 2.} Assume in addition to the example 1

$$
m(\delta) \le w + \kappa |\ln \delta|, \ 0 < \delta \le 1, w = \const > 0, \ \kappa = \const \ge 0.
$$
Then

$$
g(p) \le w + \frac{\kappa |\ln p|}{1 - p}, \ p \in (0,1),
$$
 hence

$$
P(T,u) \le w \ \sqrt{e} \ [\pi(u)]^{-\kappa} \ e^{-\phi^*(u)}, \ u > u_0.
$$

\vspace{3mm}

\ {\bf Example 3.} Let now

$$
  m(\delta) \le \delta^{-\nu}, \ \delta \in (0,1), \ \nu = \const \in (0,1);
$$
then

$$
g(p) \le p^{-\nu} \ \cdot \frac{1 - p}{1 - p^{1 - \nu}} \le  p^{-\nu} \ \cdot \frac{1}{1 - 2^{\nu - 1}}, \ p \in(0,1/2).
$$

 \ We get:

$$
P(T,u) \le C_1(\nu) \exp \left(-\phi^*(u) + C_2(\phi) \cdot [\pi(u)]^{-\nu/(\nu + 1)} \right), \ u \ge 3.
$$

 \vspace{4mm}

 \ {\bf B. Lower bounds.}\par

\vspace{3mm}

 \ The lower bounds is more simple. Namely, we will use the following fact, see
\cite{Gorskikh1},  \cite{Ostrovsky1}, pp. 25-27. For arbitrary normed:  $  \phi(0) = 0  $ Young-Orlicz function
$  \phi = \phi(\lambda), \ \lambda \in R^1 $ there exists a centered random variable $  \rho  $ having the unit norm
in the space $  \ B(\phi): \ ||\rho||B\phi = 1  $ but for which there exists a finite positive constant $  C $ and
a sequence of real numbers $  \ x(n), \ n = 1,2,3,\ldots \ $ tending to infinity as $  n \to \infty $ such that

$$
{\bf  P} (\rho > x(n)) \ge C \exp \left( - \phi^*(x(n)) \right).
$$

 \ The correspondent (continuous!) random fields $  \xi_0 = \xi_0(t), \ t \in T  $ may be constructed as a constant relative
the "time" variable  $ t: \ \xi_0(t) = \rho $ with at the same estimate for  behavior of the maximum distribution. \par

\vspace{4mm}

\section{ Exponential estimates for maximum distribution of normed sums of discontinuous random fields. }

\vspace{3mm}

 \  Let $   T = [0,1]^m  $ be again the ordinary $ m \ -  $ dimensional unit cube and let

$$
\xi = \xi(\vec{t}) = \xi( t(1), t(2), \ldots, t(m) ), \hspace{4mm} t(i) \in [0,1], \ t = \vec{t} \in T
$$
be the separable numerical valued mean zero random process (field).   Let also $ \xi_i = \xi_i(t), \ i = 1,2, \ldots $
be independent copies of the r.f. $  \xi(t). $ Denote

$$
S_n(t) = n^{-1/2} \sum_{i=1}^n \xi_i(t),\eqno(4.0)
$$
and put

$$
Q_n(T,u) =  {\bf P} ( \sup_{t \in T} S_n(t) > u), \eqno(4.1)
$$
and correspondingly

$$
Q(T,u) = \sup_n Q_n(T,u) =  \sup_n {\bf P} ( \sup_{t \in T} S_n(t) > u).\eqno(4.1a)
$$

\vspace{3mm}

 \ {\it Our goal in this section is deriving an exact non-asymptotical  exponential  decreasing to zero as $  u \to \infty  $
estimates for } $  Q_n(T,u), \ Q(T,u). $  \par

\vspace{3mm}

 \ These estimates may be used, for example, in the statistics of random processes and for building of non-asymptotical
confidence region in the uniform norm in the Monte-Carlo method for the computation of multiparameter integrals with
discontinuous integrand function, see  \cite{Grigorjeva1}. \par

\vspace{3mm}

 Let again $  \phi = \phi(\lambda) $ be Young - Orlicz  function, $  \phi(0) = \phi(0+) = 0,  $ and such that
for the  centered r.f. $ \xi(t), \ t \in [0,1]^m  $

$$
 \sup_{t \in T} {\bf E} e^{\lambda \xi(t)} \le e^{\phi(\lambda) }, \ \lambda \ge 0. \eqno(4.2)
$$

 \ It will be presumed as before the  finiteness of these function on some right-hand side non-empty neighborhood of zero. \par

 \ For instance, the function $ \phi = \phi_0(\lambda) $ may be picked by the natural way:

$$
e^{\phi_0(\lambda) } :=  \sup_{t \in T} {\bf E} e^{\lambda \xi(t)}, \ \lambda \ge 0.
$$

 \ Let us introduce the following functions

$$
\phi_n(\lambda) :=  n \ \phi(\lambda/\sqrt{n}), \hspace{5mm}  \overline{\phi}(\lambda) := \sup_n [ n \ \phi(\lambda/\sqrt{n})],
 \eqno(4.3)
$$
and correspondingly

$$
g_n(p) = g(p,\phi_n), \hspace{5mm}  \overline{g}(p) = g(p,\overline{\phi}).  \eqno(4.4)
$$

\vspace{4mm}

{\bf  Theorem B.} \par

\vspace{3mm}

\ {\bf I.} Suppose

$$
 \exists p \in (0,1)  \ \Rightarrow g_n(p) < \infty.  \eqno(4.5)
$$

\vspace{3mm}

\ Then

$$
Q_n(T,u) \le \inf_{p \in (0,1) } \exp \left[ g_n(p) - \phi_n^*(u(1-p)) \right], \ u > 0.  \eqno(4.5a)
$$

\vspace{3mm}

\ {\bf II.} Suppose

$$
 \exists p \in (0,1)  \ \Rightarrow  \overline{g}(p) < \infty.  \eqno(4.6)
$$

\vspace{3mm}

\ Then

$$
Q(T,u) \le \inf_{p \in (0,1) } \exp \left[ \overline{g}(p) - \overline{\phi}^*(u(1-p)) \right], \ u > 0.  \eqno(4.6a)
$$

\vspace{3mm}

{\bf Proof} follows immediately from theorem A, if we note by virtue of independence of r.f. $ \xi_i(\cdot)  $

$$
\sup_t {\bf E } e^{\lambda S_n(t)} = \sup_t \left[ {\bf E } e^{ \lambda \xi(t)/\sqrt{n}} \right]^n \le
e^{n \phi(\lambda /\sqrt{n})} = e^{ \phi_n(\lambda)}
$$
and consequently

$$
\sup_n \sup_t {\bf E } e^{\lambda S_n(t)} = \sup_n \sup_t \left[ {\bf E } e^{ \lambda \xi(t)/\sqrt{n}} \right]^n \le
\sup_n e^{n \phi(\lambda /\sqrt{n})} = e^{ \overline{\phi}(\lambda)}.
$$

\vspace{5mm}

 \ The convenient lower estimate for the probability $  Q(T,u)  $ one can obtain as follows:

$$
Q(T,u) \ge \max \left\{P(T,u), \ \lim_{n \to \infty} {\bf P}( \sup_t S_n(t) > u ) \right\},
$$
where the second member in curly brackets may be estimated  by means of the classical one-dimensional CLT

$$
\lim_{n \to \infty} {\bf P}( \sup_t S_n(t) > u )  \ge \lim_{n \to \infty} {\bf P}(  S_n(t_0) > u ) \ge
$$

$$
\exp(-c u^2), \ u \ge 1,
$$
where $  t_0 $ is arbitrary point inside the set $   T  $ for which $  \Var (\xi(t_0)) > 0. $ \par

 \vspace{3mm}

 \section{Concluding remarks.}

 \vspace{3mm}

 \ The exponential estimate  for distribution of absolute value for $  \xi(t) $ may be easily obtained quite
analogously by means of relation

$$
{\bf P} ( \max_{t \in T} |\xi(t)| > u ) \le {\bf P} ( \max_{t \in T} \xi(t) > u ) +
{\bf P} ( \max_{t \in T} ( -\xi(t)) > u ), \ u > 0.
$$

\end{document}